\documentclass[draft,12pt]{amsart}
\usepackage{amssymb}

\theoremstyle{plain}
\newtheorem*{theorem}{Theorem}

\theoremstyle{definition}
\newtheorem*{remark}{Remark}

\newcommand{\abs}[1]{\lvert#1\rvert}
\newcommand{\norm}[1]{\lVert#1\rVert}
\newcommand{\bigabs}[1]{\bigl\lvert#1\bigr\rvert}
\newcommand{\bignorm}[1]{\bigl\lVert#1\bigr\rVert}
\newcommand{\Bigabs}[1]{\Bigl\lvert#1\Bigr\rvert}

\renewcommand{\leq}{\leqslant}
\renewcommand{\geq}{\geqslant}

\newcommand{\term}[1]{{\textit{\textbf{#1}}}}

\def\B{B(x_0,\varepsilon)}
\def\S{S(x_0,\varepsilon)}

\DeclareMathOperator{\Range}{Range}

\begin{document}
\baselineskip 18pt

\title[Minimal vectors in Banach spaces]
      {Minimal vectors in\\
                        arbitrary Banach spaces}

\author[V.~G. Troitsky]{Vladimir~G. Troitsky} 
\address{Department of Mathematics,
         University of Alberta, Edmonton, AB, T6G\,2G1. Canada.
         }
\email{vtroitsky@math.ualberta.ca}


\begin{abstract}
  We extend the method of minimal vectors to arbitrary Banach spaces.
  It is proved, by a variant of the method, that certain
  quasinilpotent operators on arbitrary Banach spaces have
  hyperinvariant subspaces.
\end{abstract}

\maketitle

The method of {\it minimal vectors} was introduced by Ansari and Enflo
in~\cite{Ansari:98} in order to prove the existence of invariant
subspaces for certain classes of operators on a Hilbert space. Pearcy
used it in~\cite{Pearcy} to prove a version of Lomonosov's theorem.
Androulakis in~\cite{Androulakis} adapted the technique to
super-reflexive Banach spaces. In~\cite{Chalendar:A} the method was
independently generalized to reflexive Banach spaces. There has been a
hope that this technique could provide a positive solution to the
invariant subspace problem for these spaces. In this note we present a
version of the method of minimal vectors (based on~\cite{Androulakis})
that works for arbitrary Banach spaces. In particular, it applies in
the spaces where there are known examples of operators without
invariant subspaces, e.g., \cite{Enflo:76,Enflo:87,Read:84,Read:85}.
This shows that the method of minimal vectors alone cannot solve the
invariant subspace problem for ``good'' spaces.

Suppose that $X$ is a Banach space. For simplicity, we
assume that $X$ is a real Banach space, though the results can be
adapted to the complex case in the straightforward manner. In the
following, $\B$ stands for the closed ball of radius $\varepsilon$
centered at $x_0$ while $B^\circ(x_0,\varepsilon)$ stands for the open
ball, and ${\S}$ stands for the corresponding sphere.

Let $Q$ be a bounded operator on $X$. Since
we will be interested in the hyperinvariant subspaces of $Q$, we can
assume without loss of generality that $Q$ is one-to-one and has dense
range, as otherwise $\ker Q$ or $\overline{\Range Q}$ would be
hyperinvariant for $Q$. By $\{Q\}'$ we denote the commutant of $Q$.

Fix a point $x_0\neq 0$ in $X$ and a positive real
$\varepsilon<\norm{x_0}$. Let $K=Q^{-1}\B$. Clearly, $K$ is a convex
closed set. Note that $0\notin K$ and $K\neq\varnothing$ because $Q$
has dense range. Let $d=\inf_{K}\norm{z}$, then $d>0$. It is observed
in~\cite{Ansari:98,Androulakis} that if $X$ is reflexive, then there
exists $z\in K$ with $\norm{z}=d$, such a vector is called a
\term{minimal vector} for $x_0$, $\varepsilon$ and~$Q$. Even without
reflexivity condition, however, one can always find $y\in K$ with
$\norm{y}\leq 2d$, such a $y$ will be referred to as a \term{2-minimal
  vector} for $x_0$, $\varepsilon$ and $Q$.

The set $K\cap B(0,d)$
is the set of all minimal vectors, in general this set may be
empty. If $z$ is a minimal vector,
since $z\in K=Q^{-1}\B$ then $Qz\in\B$. As $z$ is an
element of minimal norm in $K$ then, in fact, $Qz\in{\S}$. 
Since $Q$ is one-to-one, we have
\begin{displaymath}
  QB(0,d)\cap\B=Q\bigl(B(0,d)\cap K)\subseteq{\S}.
\end{displaymath}
It follows that $QB(0,d)$ and $B^\circ(x_0,\varepsilon)$ are two
disjoint convex sets. Since one of them has non-empty interior, they
can be separated by a continuous linear functional (see,
e.g.,~\cite[Theorem~5.5]{Aliprantis:99}). That is, there exists a
functional $f$ with $\norm{f}=1$ and a positive real $c$
such that $f_{|QB(0,d)}\leq c$ and $f_{|B^\circ(x_0,\varepsilon)}\geq c$.
By continuity, $f_{|\B}\geq c$. We say that $f$ is a
\term{minimal functional} for $x_0$, $\varepsilon$, and $Q$.

We claim that $f(x_0)\geq\varepsilon$. Indeed, for every $x$ with
$\norm{x}\leq 1$ we have $x_0-\varepsilon x\in\B$. It follows that
$f(x_0-\varepsilon x)\geq c$, so that $f(x_0)\geq c+\varepsilon f(x)$.
Taking $\sup$ over all $x$ with $\norm{x}\leq 1$ we get $f(x_0)\geq
c+\varepsilon\norm{f}\geq\varepsilon$.

Observe that the hyperplane $Q^*f=c$ separates $K$ and $B(0,d)$.
Indeed, if $z\in B(0,d)$, then $(Q^*f)(z)=f(Qz)\leq c$, and if $z\in K$
then $Qz\in\B$ so that $(Q^*f)(z)=f(Qz)\geq c$. 
For every $z$ with $\norm{z}\leq 1$ we have $dz\in B(0,d)$, so
that $(Q^*f)(dz)\leq c$, it follows that
$\bignorm{Q^*f}\leq\frac{c}{d}$. On the other hand, for every
$\delta>0$ there exists $z\in K$ with $\norm{z}\leq d+\delta$,
then $(Q^*f)(z)\geq c\geq\frac{c}{d+\delta}\norm{z}$, whence
$\bignorm{Q^*f}\geq\frac{c}{d+\delta}$. It follows that
$\bignorm{Q^*f}=\frac{c}{d}$. For every $z\in K$ we have
$(Q^*f)(z)\geq c=d\bignorm{Q^*f}$. In particular, if $y$ is
a 2-minimal vector then
\begin{equation}
  \label{eq:norm-att}
  (Q^*f)(y)\geq\tfrac{1}{2}\bignorm{Q^*f}\norm{y}.
\end{equation}

We proceed to the main theorem

\begin{theorem}\label{t:main}
  Let $Q$ be a quasinilpotent operator on a Banach space $X$, and
  suppose that there exists a closed ball $B$ such that $0\notin B$ and
  for every sequence $(x_n)$ in $B$ there is a subsequence $(x_{n_i})$
  and a sequence $(K_i)$ in $\{Q\}'$ such that $\norm{K_i}\leq 1$ and
  $(K_ix_{n_i})$ converges in norm to a non-zero vector. Then $Q$ has
  a hyperinvariant subspace.
\end{theorem}

\begin{remark}
  The hypothesis of the theorem is slightly weaker than the condition ($*$)
  in~\cite{Androulakis}, where it is required that for every
  $\varepsilon\in(0,1)$ there exists $x_0$ of norm one such that the
  ball $\B$ satisfies the rest of the condition.
\end{remark}

\begin{proof}
  Without loss of generality $Q$ is one-to-one and has dense range.
  Let $x_0\neq 0$ and $\varepsilon\in\bigl(0,\norm{x_0}\bigr)$ be such
  that $B=\B$.  For every $n\geq 1$ choose a 2-minimal vector $y_n$
  and a minimal functional $f_n$ for $x_0$, $\varepsilon$, and $Q^n$.
  
  Since $Q$ is quasinilpotent, there is a subsequence $(y_{n_i})$ such
  that $\frac{\norm{y_{n_i-1}}}{\norm{y_{n_i}}}\to 0$. Indeed,
  otherwise there would exist $\delta>0$ such that
  $\frac{\norm{y_{n-1}}}{\norm{y_n}}>\delta$ for all $n$, so
  that
  $\norm{y_1}\geq\delta\norm{y_2}\geq\dots\geq\delta^n\norm{y_{n+1}}$.
  Since $Q^ny_{n+1}\in Q^{-1}B$ then
  $$\bignorm{Q^ny_{n+1}}\geq d\geq\frac{\norm{y_1}}{2}\geq
  \frac{\delta^n}{2}\norm{y_{n+1}}.$$
  It follows that
  $\norm{Q^n}\geq\delta^n/2$, which contradicts the quasinilpotence
  of~$Q$.
  
  Since $\norm{f_{n_i}}=1$ for all $i$, we can assume (by passing to
  a further subsequence), that $(f_{n_i})$ weak*-converges to some
  $g\in X^*$. Since $f_{n}(x_0)\geq\varepsilon$ for all $n$, it follows
  that $g(x_0)\geq\varepsilon$. In particular, $g\neq 0$.
  
  Consider the sequence $(Q^{n_i-1}y_{n_i-1})_{i=1}^\infty$. It is
  contained in $B$, so that
  by passing to yet a further subsequence, if necessary, we find 
  a sequence $(K_i)$ in $\{Q\}'$ such that $\norm{K_i}\leq
  1$ and $K_iQ^{n_i-1}y_{n_i-1}$ converges in norm to some $w\neq 0$.
  Put
  $$Y=\{Q\}'Qw=\bigl\{TQw\mid T\in\{Q\}'\bigr\}.$$
  One can easily
  verify that $Y$ is a linear subspace of $X$ invariant under
  $\{Q\}'$. Notice that $Y$ is non-trivial because $Q$ is one-to-one
  and $0\neq Qw\in Y$. We will show that $Y\subseteq\ker g$, so that
  $\overline{Y}$ is a proper $Q$-hyper\-in\-variant subspace.
  
  Take $T\in\{Q\}'$; we will show that $g(TQw)=0$. It follows from
  (\ref{eq:norm-att}) that $\bigl(Q^{*n_i}f_{n_i}\bigr)(y_{n_i})\neq
  0$ for every $i$, so that $X={\rm span}\{y_{n_i}\}\oplus
  \ker\bigl(Q^{*n_i}f_{n_i}\bigr)$.  Then one can write
  $TK_iy_{n_i-1}=\alpha_iy_{n_i}+r_i$, where $\alpha_i$ is a scalar
  and $r_i\in \ker\bigr(Q^{*n_i}f_{n_i}\bigl)$. We claim that
  $\alpha_i\to 0$.  Indeed,
  \begin{equation}\label{eq:reduce}
    \bigl(Q^{*n_i}f_{n_i}\bigr)\bigl(TK_iy_{n_i-1}\bigr)=
    \alpha_i\bigl(Q^{*n_i}f_{n_i}\bigr)(y_{n_i}),
  \end{equation}
  and, combining this with (\ref{eq:norm-att}) we get
  \begin{equation}\label{eq:abs1}
    \bigabs{\bigl(Q^{*n_i}f_{n_i}\bigr)\bigl(TK_iy_{n_i-1}\bigr)}\geq
    \frac{\abs{\alpha_i}}{2}\bignorm{Q^{*n_i}f_{n_i}}\norm{y_{n_i}}.
  \end{equation}
  On the other hand,
  \begin{equation}\label{eq:abs2}
    \bigabs{\bigl(Q^{*n_i}f_{n_i}\bigr)\bigl(TK_iy_{n_i-1}\bigr)}\leq
    \bignorm{Q^{*n_i}f_{n_i}}\cdot\norm{T}\cdot\norm{y_{n_i-1}}.
  \end{equation}
  It follows from (\ref{eq:abs1}) and (\ref{eq:abs2}) that
  $$\abs{\alpha_i}\leq2\norm{T}\frac{\norm{y_{n_i-1}}}{\norm{y_{n_i}}}\to
    0.$$
  Then (\ref{eq:reduce}) yields that
  \begin{multline*}
    \Bigabs{f_{n_i}\bigl(Q^{n_i}TK_iy_{n_i-1}\bigr)}=
    \Bigabs{\alpha_if_{n_i}\bigl(Q^{n_i}y_{n_i}\bigr)}\\\leq
    \abs{\alpha_i}\cdot\norm{f_{n_i}}\cdot\norm{Q^{n_i}y_{n_i}}\leq
    \abs{\alpha_i}\cdot 1\cdot\bigl(\norm{x_0}+\varepsilon\bigr)\to 0,
  \end{multline*}
  so that $f_{n_i}\bigl(Q^{n_i}TK_iy_{n_i-1}\bigr)\to 0$.  On the
  other hand, since $T,K_i\in\{Q\}'$ we have
  $$Q^{n_i}TK_iy_{n_i-1}=TQK_iQ^{n_i-1}y_{n_i-1}\to TQw$$ in norm, while
  $f_{n_i}\xrightarrow{w^*}g$, so that
  $f_{n_i}\bigl(Q^{n_i}TK_iy_{n_i-1}\bigr)\to g(TQw)$. Hence, $g(TQw)=0$.
\end{proof}

\bigskip

Clearly, the argument will work as well for $\lambda$-minimal vectors
for any $\lambda>1$.

Suppose that $Q$ is a quasinilpotent operator commuting with a
compact operator $K$. Then $Q$ satisfies the hypothesis of
Theorem~\ref{t:main}. Indeed, without loss of generality
$\norm{K}=1$. Fix $\varepsilon=\frac{1}{3}$, there exists $x_0$
with $\norm{x_0}=1$ such that $\norm{Kx_0}\geq\frac{2}{3}$, then
$0\notin K\B$. For every sequence $(x_n)$ in $\B$, the sequence
$(Kx_n)$ has a convergent subsequence $(Kx_{n_i})$. Take $K_i=K$ for all
$i$; since $0\notin K\B$ then $\lim_iK_ix_{n_i}\neq 0$.
It follows from Theorem~\ref{t:main} that {\it if $Q$ is a
  quasinilpotent operator on a real or complex Banach space commuting
  with a non-zero compact operator, then $Q$ has a hyperinvariant
  subspace.} This fact is not new though: for complex Banach spaces it is
a special case of the celebrated Lomonosov's
theorem~\cite{Lomonosov:73}, and for real Banach spaces it follows
from Theorem~2 of~\cite{Hooker:81}.

%To apply Theorem~2 of~\cite{Hooker:81} we need to show that $Q$ does
%not satisfy any irreducible polynomial relation. Indeed, suppose that
%$p(Q)=0$ for some polynomial $p(t)$. Decompose $p(t)$ into a product
%$p(t)=p_1(t)p_2(t)\dots p_n(t)$ where the factors have degree at most
%2. Show that $n=1$. Indeed, if $n>1$ then since $p(Q)$ is irreducible,
%we have $p_n(T)\neq 0$. But since $\Range p_n(T)$ is
%$Q$-hyperinvariant, it is dense. Then $p_1(T)p_2(T)\dots p_{n-1}(T)$
%vanishes on a dense set, hence it is zero, contradiction.  Thus,
%$n=1$. But $p(t)$ cannot be of degree one because then $Q$ would have
%an eigenvector. Hence, $p(t)$ is a quadratic with roots
%$\lambda_1,\lambda_2\in\mathbb C\setminus\mathbb R$. Then the
%complixification $Q_{\mathbb C}$ satisfies $(\lambda_1I-Q_{\mathbb
%  C})(\lambda_2I-Q_{\mathbb C})=0$. Pick any non-zero $x\in X$.  If
%$(\lambda_2I-Q_{\mathbb C})x=0$ then $\lambda_2$ is an eigenvalue of
%$Q_{\mathbb C}$. If not, then $(\lambda_1I-Q_{\mathbb
%  C})$ vanishes on $(\lambda_2I-Q_{\mathbb C})x\neq 0$, hence
%$\lambda_1$ is an eigenvalue of $Q_{\mathbb C}$. In either case,
%$\sigma(Q_{\mathbb C})$ contains a non-zero element. But
%$\sigma(Q_{\mathbb C})=\sigma(Q)$, which contradicts the
%quiasinilpotence of $Q$.

\textbf{Acknowledgements.}  Thanks are due to A.~Litvak,
N.~Tomczak-Jaegermann, and R.~Vershynin for their interest in the work
and helpful suggestions. The author would also like to thank the
Department of Mathematics of the University of Alberta for its support
and hospitality.

\nocite{Radjavi:73}

\end{document}